\RequirePackage{fix-cm}
\documentclass[smallextended]{svjour3}       % onecolumn (second format)
\smartqed  % flush right qed marks, e.g. at end of proof
\usepackage{graphicx}
\usepackage{longtable}
\journalname{SEMR}
\begin{document}

\title{Mutually touching infinite cylinders in the 3D world of lines%\thanks{Grants or other notes
}

\titlerunning{Mutually touching cylinders}        % if too long for running head

\author{Peter V Pikhitsa         \and
       Stanislaw Pikhitsa  %etc.
}


\institute{Peter V. Pikhitsa \at
              Seoul National University, 151-744 Seoul, Korea \\
              Tel.: +822-88-01694\\
              Fax: +822-88-06671\\
              \email{peter@snu.ac.kr}           %  \\
           \and
           Stanislaw Pikhitsa \at
              Odessa National University, 65082 Odessa, Ukraine
}

\date{Received: date / Accepted: date}

\maketitle

\begin{abstract}
Recently we gave arguments that only two unique topologically different configurations of 7 equal all mutually touching round cylinders (the configurations being mirror reflections of each other) are possible in 3D, although a whole world of configurations is possible already for round cylinders of arbitrary radii. It was found that as many as 9 round cylinders (all mutually touching) are possible in 3D while the upper bound for arbitrary cylinders was estimated to be not more than 14 under plausible arguments. Now by using the chirality and Ring matrices that we introduced earlier for the topological classification of line configurations, we have given arguments that the maximal number of mutually touching straight infinite cylinders of arbitrary cross-section (provided that its boundary is a smooth curve) in 3D cannot exceed 10. We generated numerically several configurations of 10 cylinders, restricting ourselves with elliptic cylinders.  Configurations of 8 and 9 equal elliptic cylinders (all in mutually touching) are generated numerically as well.  A possibility and restriction of continuous transformations from elliptic into round cylinder configurations are discussed. Some curious results concerning the properties of the chiral matrix (which coincides with Seidel's adjacency matrix important for the Graph theory) are presented.

\keywords{mutual touching \and infinite cylinders \and ultimate configurations \and topology}
\subclass{52C17 \and 65H04 \and 57M99}
\end{abstract}

\section{\textbf{Introduction}}
\label{intro}
Entanglement in 3D plays an important role in the theory of knots and their invariants. Nowadays the knot invariants are used in such fundamental areas as the field theories and string theories.  Yet, seemingly standing aloof, topology and geometry of lines have not been developed as much probably because up to recent times there has been little evidence for any practical use.  In spite of this, the entanglement of lines in 3D is far from trivial and certainly deserves studying. Especially interesting are "coated" straight lines --- the cylinders of arbitrary cross-section provided that its boundary is a smooth curve.

As a prominent challenge there has been the problem of seven mutually contacting equal round cylinders first solved numerically in \cite{Ref1} and called 7*-knot in \cite{Ref2}, \cite{Ref3}. The 7*-knot was rediscovered five years later independently in \cite{Ref4}. We have called the configuration of $n$ mutually touching cylinders an $n$-knot \cite{Ref5}, bearing in mind the entanglement but other proper names are possible. The closest to the sense is the $n$-cross, the name that could combine both the old name for a rod and the tree-like structure of the configuration. In the present paper we are still using the name of $n$-knot only for the mutually touching cylinders. The name $n$-cross we reserve for a configuration of $n$ cylinders that are not obliged to be in mutually touching. Of course, $n$-crosses include $n$-knots as subsets.

\begin{figure}[h]
  \includegraphics [width=0.85\textwidth]{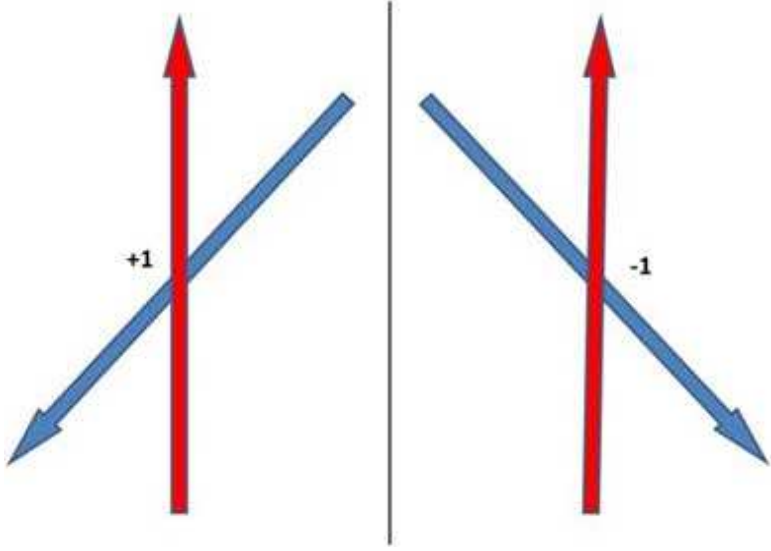}
\caption{ The chirality +1 (left) and -1 (right) of the configuration of two oriented lines in 3D. }
\label{fig:1}       % Give a unique label
\end{figure}

Another challenge could be a crystal of interwoven cylinders (we will call this broader structure that includes $n$-knots and $n$-crosses a lineation in analogy with tessellation) as a kind of extension of liquid crystal paradigm but with an intricate entangled topology and possibly auxetic mechanical properties \cite{Ref5}, \cite{Ref1}. At last, the knowledge of line configuration perhaps could be useful for characterization of so-called stick knots. Other prospects are to be discovered.
Lineation, $n$-cross and $n$-knot need an adequate formulation to unambiguously distinguish line configurations except maybe the degenerate ones. This formulation was started being created in \cite{Ref2}, \cite{Ref3}. Here we briefly introduce it to the reader. For $n$ oriented lines ($n$-cross) one needs a normalized chirality matrix $P_{i,k}$ \cite{Ref2} - a symmetric matrix with off-diagonal entries +1 and -1, and 0 on the diagonal. For any given $n$-cross one can easily construct the matrix while inspecting a configuration by following the simple rule depicted in Fig. \ref{fig:1}: when the rotation from the line that covers the other goes counterclockwise then the entry is +1, otherwise -1.

The determinant of $P_{i,k}$ is a topological invariant of an $n$-cross or of an $n$-knot. The most important for characterization of $n$-knots are two chirality matrices:

 \begin{equation}
K5=\left( \begin{array}{rrrrr}
0&1&1&1&1\\
1&0&1&1&1\\
1&1&0&1&1\\
1&1&1&0&1\\
1&1&1&1&0 \end{array} \right)
\label{eq:1}
\end{equation}

(well-known from the Graph Theory as one of Kuratowski's matrices) and the one with the determinant 250 (or -250 for the mirror one)

 \begin{equation}
P250=\left( \begin{array}{rrrrrrr}
0&1&1&-1&-1&-1&1\\
1&0&1&1&1&-1&1\\
1&1&0&1&1&-1&1\\
-1&1&1&0&1&-1&-1\\
-1&1&1&1&0&1&1\\
-1&-1&-1&-1&1&0&1\\
1&1&1&-1&1&1&0 \end{array} \right).
\label{eq:2}
\end{equation}

Both matrices in Eq. (\ref{eq:1}) and Eq. (\ref{eq:2}) correspond to configurations impossible for mutually touching as was proved in \cite{Ref3}. If one finds Eq. (\ref{eq:1}) or Eq. (\ref{eq:2}) or their equivalent as the sub-matrices in any chirality matrix, then the total configuration is impossible for mutually touching.

Note aside from the main discussion that there are interesting properties of the "extreme" matrices with the maximum possible determinant which is 250 for rank 7.  For the matrix of Eq. (\ref{eq:2}) we found that the determinant $|P_{i,k} (r_i+r_k)|$ can turn into zero for all non-negative $r_i$ only when all $r_i=0$ \cite{Ref3}. It was proved analytically in direct calculation by finding the series of all positive terms for the determinant. We claim that a general property of having the same sign of the determinant for all $r_i > 0$   takes place for extreme chirality matrices of any rank. We proved it analytically for ranks up to 9 with Mathcad11 and we hope to give a general proof elsewhere by using the method of \cite{Ref6} and \cite{Ref7} which was performed on the matrix with $P_{i,k}=1$ and $P_{i,i}=0$ to find $|P_{i,k} (r_i+r_k)|$ analytically.

When the rank of chirality matrices grows there is a crowding of forbidden $K5$ and $P250$ sub-matrices that both reduce the maximum rank of the allowed matrix of chirality down to 14 \cite{Ref3}. That means that not more than 14 cylinders could be driven into mutually touching. Without $P250$ the maximum rank of the matrix that does not contain $K5$ would be 18 \cite{Ref3}.
It is interesting that both allowed matrices of rank 14 and 18 are extreme in a sense that they belong to the series of Seidel's matrices with the largest determinants and the simplest characteristic polynomial. We put forward a conjecture that for $n=0,1,2 \ldots $ the extreme matrices $P$ of the rank $N=4n+2$ have the characteristic polynomial $(x^2-(N-1))^{N/2}$ which for $n=3$ and $n=4$ coincide with the previously reported $(x^2-13)^7$ and $(x^2-17)^9$, respectively \cite{Ref3}.
The simple polynomial was found by noticing that for the extreme matrix of the rank $N$ defined above there is an equation $P^2=(N-1)I$, where $I$ is the identity matrix of rank $N$. We have not proved the general formula  but we plan to do it with the methods of \cite{Ref6},\cite{Ref7}. For $n=1$ one gets $(x^2-5)^3$ so that the determinant of the extreme matrix of rank 6 is $-125$. As we show below, this matrix plays a crucial role in the analysis of the $n$-crosses forbidden for mutually touching of the cylinders.
 Not normalized chirality matrix is defined as \cite{Ref2}

 \begin{equation}
{\cal P}_{i,k}\equiv \left( \begin{array}{cccc}
0&R_{ik} [{\textbf{n}_i} \times {\textbf{n}_k}]&\cdots&\cdots\\
R_{ik} [{\textbf{n}_i} \times {\textbf{n}_k}]&0&\cdots&\vdots\\
\vdots&\vdots&\ddots&\vdots\\
\vdots&\cdots&\cdots&0 \end{array} \right)
\label{eq:3}
\end{equation}

where $\textbf{n}_i$ is a unit vector along the $i$th line and $R_{ik}$ is the vector of the shortest distance between the  $i$th and the $k$th lines, collinear with $[{\textbf{n}_i} \times {\textbf{n}_k}]$. This matrix coincides with the matrix of products of the line moments introduced for seven lines in \cite{Ref7} (see item 17 on page 194 therein). There is an interesting property of the determinant of such a matrix: the determinant is always zero for seven and more lines. It happens because a line can be characterized by a vector of six values: three components of the unit vector of the line direction and three components of the line moment with respect to some point. With the line moment being orthogonal to the unit vector, the matrix of Eq. (\ref{eq:3}) is formed as a product of a matrix $(6 \times n)$ to the transposed one $(n \times 6)$. The resulting $(n \times n)$ matrix inevitably has zero determinant when $n>6$. Being the Gram determinant that gives the square of the volume of the parallelotope of the six dimensional vectors it is naturally zero when the number of the vectors is more than 6.  This brings an interesting separation of the configurations of the round cylinders from a broader set of configurations of arbitrary shaped cylinders which we discuss below. Finally, the normalized chirality matrix $P_{i,k}$ discussed above is constructed from Eq. (\ref{eq:3}) by just taking the sign of each entry $P_{i,k}=sign({\cal P}_{i,k})$.

\begin{figure}[h]
  \includegraphics [width=0.85\textwidth]{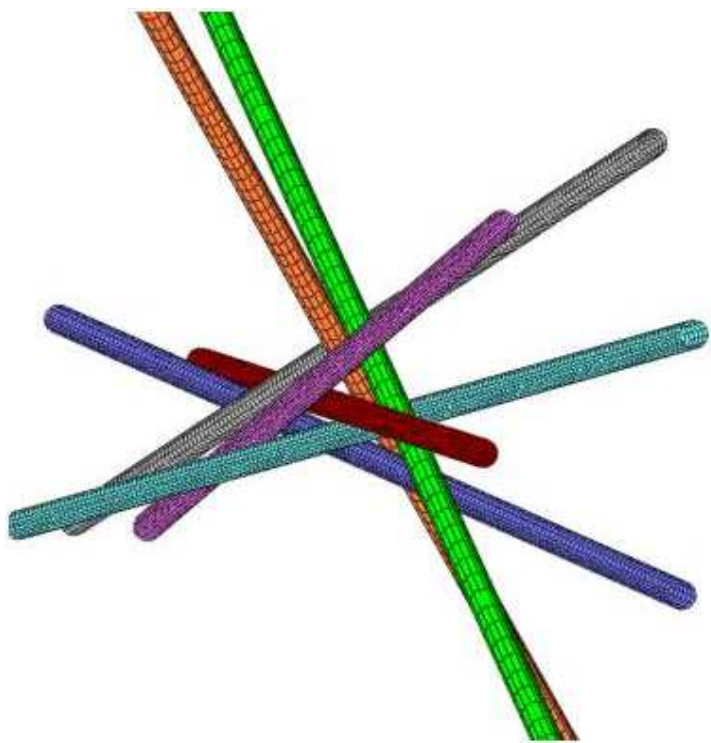}
\caption{The 7*-knot of all equal mutually pairwise touching round cylinders.}
\label{fig:2}       % Give a unique label
\end{figure}

What remains to discuss in the Introduction part is the Ring matrix $R_{i,k}$ that completes the topological characterization of a lineation, $n$-cross and $n$-knot along with the chirality matrix. Although it has a proper rule of calculation \cite{Ref3}, its powerful feature is the possibility to construct the Ring matrix by a simple inspection of a given configuration of an $n$-cross. To make this one has to watch the configuration along a chosen cylinder, say, of number 0 which axis is projected as a point. Axes of all other $n-1$ cylinders are projected on the plane of view as $n-1$ lines. Their distribution around the point is what forms the row of the Ring matrix. Let us take the first line and count the number of triangles which have the first line as a lateral side and which have the point 0 in their inside so that the line 0 is said to be encaged by the other lines. This number makes the entry  $R_{0,1}$ . Proceeding in the same manner one obtains the full row $R_{0,i}$, where $i=1,\ldots n-1$. The diagonal $R_{0,0}$   is taken to be 0. Then we watch the configuration along the axis of the cylinder number 1, etc. Eventually, the Ring matrix will have null diagonal and non-zero entries when the lines are encaged/entangled. Otherwise the rows remain all zeroes and such lines are called free. An example of the Ring matrix for the unique configuration (7*-knot, Fig. \ref{fig:2}) of seven equal round cylinders (all mutually touching) was given in \cite{Ref3} along with the chirality matrix:

 \begin{equation}
R=\left( \begin{array}{ccccccc}
0&1&1&4&1&1&4\\
4&0&4&4&4&4&4\\
0&0&0&0&0&0&0\\
0&0&0&0&0&0&0\\
1&1&4&4&0&1&1\\
1&1&1&1&4&0&4\\
0&0&0&0&0&0&0 \end{array} \right)
\label{eq:4}
\end{equation}

\begin{equation}
P=\left( \begin{array}{rrrrrrr}
0&1&1&1&1&1&1\\
1&0&1&1&1&-1&1\\
1&1&0&-1&-1&-1&1\\
1&1&-1&0&-1&1&1\\
1&1&-1&-1&0&1&-1\\
1&-1&-1&1&1&0&1\\
1&1&1&1&-1&1&0 \end{array} \right)
\label{eq:5}
\end{equation}

It is interesting to notice that the sum of the entries in a row in the Ring matrix divided by 3 gives the number of rings (triangles) that enclose the line with the same number as the number of the row. For the example given in Eq. (\ref{eq:5}) one can check that the row with all 4s gives 8 rings around the red cylinder which is easy to check by observing the image in Fig. \ref{fig:2}. Also note that the maximum number of rings that encage a line in an $n$-cross is given by the formulae:
$(n-2)(n-1)n/24$  when $n$ is even;
$(n-3)(n^2-1)/24$  when $n$ is odd.

Finally, we suggested a numerical invariant \cite{Ref3} that satisfactorily describes configurations along with the determinant of the chirality matrix. Certainly, one can suggest many kinds of invariants which are more or less devoid of singularities or repetitions for different configurations. This very invariant is
\begin{equation}
\wp(P,R)=tr[Q(I-R)^{-1}]
\label{eq:6}
\end{equation}
where $I$  is the identity matrix and we introduced a new matrix $Q$ instead of $P$ because the chirality matrix is sensitive to the directions of the directed lines. Matrix  $Q$ is invariant with respect to the line direction and formed from P as follows (we are using Eq. (\ref{eq:5}) as an example):

One has to transform the zero row in $P$ into all +1s except the diagonal term $P_{0,0}$ (in the example given it is already so; also mind that we calculate the index of row/column starting from zero) to make a matrix $A_0$ and to sum up numbers in the corresponding columns and put the sums into the zero row. Then one has to transform the first row in $P$ into all +1s except the diagonal term (here by reverting sign in the fifth column and then in the fifth row) to make a matrix $A_1$ then to sum up numbers in the corresponding columns and put the sums into the first row, etc.

After proceeding in the same way through the whole matrix of Eq. (\ref{eq:5}) one gets the symmetric matrix:
\begin{equation}
Q=\left( \begin{array}{rrrrrrr}
6&4&0&2&0&2&4\\
4&6&2&0&-2&-2&2\\
0&2&6&-2&0&2&2\\
2&0&-2&6&-2&2&4\\
0&-2&0&-2&6&0&0\\
2&-2&2&2&0&6&0\\
4&2&2&4&0&0&6 \end{array} \right)
\label{eq:7}
\end{equation}

so that the invariant can be calculated from Eq. (\ref{eq:6}) with Eqs. (\ref{eq:4}) and (\ref{eq:7}) to be $\wp(P,R)=14.07317\ldots$ and $\wp(-P,R)=7.80487\ldots$ for the mirror configuration which is characterized by the negative sign before the chirality matrix with the determinant -18 \cite{Ref2}, \cite{Ref3}. Being equipped with both the determinant and invariant one can characterize any possible $n$-knot or $n$-cross as it was done in \cite{Ref3}.  In the present paper we use the invariant in its full capacity as a research tool for finding specific configurations needed.

Note that instead of matrix $Q$ one could use another matrix which is the matrix of spirality $S$:
$S_{i,j}=-P_{i,j}  (\textbf{n}_i \cdot \textbf{n}_j)/|\textbf{n}_i \cdot \textbf{n}_j | $.
 It can substitute $Q$ in Eq. (\ref{eq:6}) and produce an invariant that becomes farther from pure topology and closer to geometry (as orthogonality of lines plays a role) and distinguishes more details of configurations, especially for mirror configurations. Spirality is independent of the orientation of lines. In the left panel in Fig. \ref{fig:1} it is +1 because watching from the angle between lines less than the right angle, non-oriented lines look rotating counterclockwise, though the non-oriented lines in the right panel look rotating clockwise and their spirality is -1. Mirror reflection changes the sign.

Yet, in the present paper we will use only the chirality for quantitative description thus staying on firm topological grounds. Still we additionally introduce here a new invariant with stronger discriminating capabilities for entangled lines. It is based on replacement of the matrix $Q(P)$ by the matrix $Qn(P,R)$, which interweaves the Ring matrix $R$ with the chirality matrix $P$, and by using $I/2$ instead of the identity matrix $I$ in Eq. (\ref{eq:6}) to avoid adverse coincidences of numbers. The interweaving with $R$ is made by putting weights to matrices $A_i$ involved in constructing $Q$ as described above and can be formalized as
 \begin{equation}
\wp n(P,R)=tr[Qn(I/2-R)^{-1}]
\label{eq:6a}
\end{equation}
 where
$Qn=3^{-6} {\sum_{i=0}}^{rank(R)-1}Rv_i  Rh_i  A_i$,     $Rv_i={\sum_{j=0}}^{rank(R)-1} R_{i,j}$ ,  $Rh_i ={\sum_{j=0}}^{rank(R)-1}R_{j,i} $.
It is important that $Qn$ preserves the same property as the matrix $Q$  \cite{Ref3}: $Qn(P,R)+Qn(-P,R)=const$ independent of $P$. This means that owing to linearity of Eq. (\ref{eq:6}) and Eq. (\ref{eq:6a}) with respect to $Q$ and $Qn$ the sum of direct and mirror invariants for any configuration with different chirality matrices but with the same Ring matrix is the same. We will illustrate this property below on important examples.

\section{\textbf{The bottleneck chirality matrix for $n$-knots at $n>10$ }  }
\label{sec:1}

As we mentioned above, sorting out chirality matrices with $K5$ and $P250$ led to the reduction of possible chirality matrices for $n>10$ to one single matrix of rank 11 (except the mirror one), one of rank 12, two of rank 13, and finally, one of rank 14 \cite{Ref3}. All $n$-knots with $n>14$ are impossible.

Note that in the present work we recalculated all relevant matrices using a direct matrix comparison procedure, while in \cite{Ref3} we used a simplified procedure where matrices were compared according to their eigenvalues. The procedure with eigenvalues may be deficient and needs more elaboration as it was shown in \cite{Ref8}. However, after recalculating with now direct procedure we confirmed our previous findings of \cite{Ref3}.  In Appendix 1 we give an analytical proof that any chirality matrix of rank more than 18 contains $K5$ and thus any $n$-knot with $n>18$ is forbidden. The proof uses the Ramsey's number $R(4,4)=18$ and we thank the anonymous referee for pointing to this number.

Meantime, we have found several of the existing 10-knot configurations (see below). However, we have been unable to construct any 11-knot out of 11-crosses and suspected that it is impossible. Then all other matrices of larger rank belong to forbidden configurations because they all contain the one with the rank 11 as a sub-matrix. In its turn, the chirality matrix of rank 11 that survived after sorting out $K5$ and $P250$ may also contain some forbidden sub-matrix that we will have to find. If so, then it will prove that only 10 arbitrary cylinders could be driven into mutually touching in 3D.

At the first glance the task looks immensely intricate and complicated. Below we give a series of heuristic arguments that demonstrate that 11-knot may not exist, and which by no means are rigorous given the complexity of the advancement through the unexplored area. We have had an experience with $P250$ matrix of rank 7 and even at such a low rank it could have been an impossible task to analyze the configuration to prove that it is forbidden for mutually touching, had not there been a specific symmetry of the matrix and the corresponding 7-cross configuration which we called "Equal Environment" where two cylinders were in topologically equivalent position (see the details in \cite{Ref3}).

Additionally, it contained two sub-matrices of rank 6 with the determinant $-125$ being the extreme matrix as discussed above. Then the presence of several sub-matrices of determinant $-125$ may be strong indication for the forbiddance of the configuration.   Therefore we carefully inspected the matrix of rank 11 given in Eq. (\ref{eq:8}) as to the presence of suspicious sub-matrices that may be forbidden, first of all the ones containing several matrices of rank 6 with the determinant $-125$.

 \begin{equation}
M11=\left( \begin{array}{rrrrrrrrrrr}
0&1&-1&1&1&-1&1&-1&1&-1&-1\\
1&0&1&1&-1&1&-1&1&1&-1&-1\\
-1&1&0&1&1&1&-1&1&-1&-1&-1\\
1&1&1&0&1&1&1&-1&-1&-1&1\\
1&-1&1&1&0&1&1&-1&-1&1&-1\\
-1&1&1&1&1&0&1&1&1&1&1\\
1&-1&-1&1&1&1&0&1&1&1&1\\
-1&1&1&-1&-1&1&1&0 &1&1&-1\\
1&1&-1&-1&-1&1&1&1 &0&-1&1\\
-1&-1&-1&-1&1&1&1&1 &-1&0&-1\\
-1&-1&-1&1&-1&1&1&-1 &1&-1&0\end{array} \right)
\label{eq:8}
\end{equation}

$M11$ has the determinant $57122=2\times13^4$ (or its negative for the mirror one) and the characteristic polynomial $-(x-2)(x+1)^2 (x^2-13)^4$.  We paid attention to a diagonal sub-matrix of $M11$ of rank 8 with determinant $1625=13\times5^3$ and it paid. $P1625$, given in Eq. (\ref{eq:9}), is the only matrix of rank 8 that contains four matrices of rank 6 with the determinant $-125$ and does not contain $K5$ and $P250$.
\begin{equation}
P1625=\left( \begin{array}{rrrrrrrr}
0&-1&-1&-1&-1&1&-1&-1\\
-1&0&-1&1&1&1&-1&1\\
-1&-1&0&-1&-1&-1&-1&1\\
-1&1&-1&0&-1&1&1&1\\
-1&1&-1&-1&0&-1&-1&-1\\
1&1&-1&1&-1&0&-1&1\\
-1&-1&-1&1&-1&-1&0&-1\\
-1&1&1&1&-1&1&-1&0 \end{array} \right)
\label{eq:9}
\end{equation}
To see the four matrices explicitly one has to rearrange the direction of lines and make some permutations between them. Such transformations do not change the determinant but show a hidden symmetry of $P1625$:
\begin{equation}
P1625=\left( \begin{array}{rrrrrrrr}
0&-1&-1&-1&-1&-1&1&-1\\
-1&0&-1&-1&-1&-1&-1&1\\
-1&-1&0&-1&1&-1&-1&-1\\
-1&-1&-1&0&-1&1&-1&-1\\
-1&-1&1&-1&0&1&1&1\\
-1&-1&-1&1&1&0&1&1\\
1&-1&-1&-1&1&1&0&1\\
-1&1&-1&-1&1&1&1&0 \end{array} \right)
\label{eq:10}
\end{equation}

In Eq. (\ref{eq:10}) one can notice two distinct diagonal blocks $4\times4$ filled with -1 or 1, correspondingly. Two identical off-diagonal $4\times4$ blocks beside -1 contain a single 1 in each row so that 1 goes in a rotation cycle from position $P_{0,6}$ to $P_{1,7}$, then to $P_{2,4}$, then to $P_{3,5}$ that implied the rotation symmetry $C_4$ of two cycles of cylinders with numbers 0,1,2,3 and 6,7,4,5 in a sandwich configuration punctured by the $C_4$  axis.  If one removes, for instance, cylinders 0 and 6, the remained configuration will have the matrix of the determinant $-125$ that we marked as $Pm125$:
\begin{equation}
Pm125=\left( \begin{array}{rrrrrr}
0&-1&-1&-1&-1&1\\
-1&0&-1&1&-1&-1\\
-1&-1&0&-1&1&-1\\
-1&1&-1&0&1&1\\
-1&-1&1&1&0&1\\
1&-1&-1&1&1&0 \end{array} \right)
\label{eq:11}
\end{equation}

It is understandable that there are exactly four such sub-matrices with the determinant $-125$ for each pair of cylinders chosen in Eq. (\ref{eq:10}). Now Eq. (\ref{eq:11}) implies the rotation symmetry $C_3$  in a sandwich configuration of two rings of three cylinders, punctured by $C_3$  axis.

\section{\textbf{6-crosses and 6-knots with determinant $-125$}}
\label{sec:3}
We accumulated statistically all possible configurations of the 6-cross with the determinant $-125$. Specifically, as configurations we generated randomly oriented lines inside a restricted cubic volume, the lines being the axes of the cylinders. Out of the configurations we selected those with the determinant $-125$ and calculated their invariants. The parameters of configurations were also tried as initial conditions for equations of mutually contacting cylinders in attempts to find 6-knots.  The list of their invariants is given in Appendix 2. The list was generated by successive comparison of a configuration newly found with those previously found and by accumulation of the data obtained. After the list of invariants stopped growing we interrupted the procedure after some reasonable delay. We also studied the 6-cross configurations with exact $C_3$ rotational symmetry and found 7 distinct configurations given in Appendix 3. Here the "short" list of invariants is reliably exhaustive. The invariants for mirror configurations coincide with the original ones. This is an accidental property of the chirality matrix with the determinant $-125$: the invariant introduced in \cite{Ref3} and discussed in the Introduction can only discriminate configurations by their Ring matrix in this case because the matrix $Q$ is trivial:
 $$ Q=\left( \begin{array}{cccccc}
5&1&1&1&1&1\\
1&5&1&1&1&1\\
1&1&5&1&1&1\\
1&1&1&5&1&1\\
1&1&1&1&5&1\\
1&1&1&1&1&5 \end{array} \right).$$

However, after a certain modification which removes the degeneracy (for example, if one uses the matrix $Qn$ introduced above and given in Appendices 2 and 3), the invariant could be different for mirror case. Here we will also show directly the difference between important configurations. After considerable numerical calculations (it is believed, though not proved, that the invariants thus obtained exhaust all possible invariants and configurations for the determinant $-125$), we noticed that of all calculated 6-crosses only two configurations, both with the identical invariant $9.6666\ldots$ can produce 6-knots. These specific configurations (which belong to the most probable ones as seen from the list in Appendix 2) may have realizations with the exact $C_3$  rotational symmetry. Their mirror invariant is also $9.6666\ldots$. The Ring matrix of these configurations is given below:
  \begin{equation}
R=\left( \begin{array}{cccccc}
0&1&1&3&3&1\\
0&0&0&0&0&0\\
3&1&0&1&1&3\\
0&0&0&0&0&0\\
1&3&3&1&0&1\\
0&0&0&0&0&0 \end{array} \right)
\label{eq:12}
\end{equation}
and it is one and the same for both configurations. It is clearly seen from Appendix 2 where the sum of new invariants for the direct and mirror configurations is given in the sixth column. One can see that it takes the same value $-1.06621$ for both $9.6666\ldots$ configurations.

\begin{figure}[h]
  \includegraphics [width=0.85\textwidth]{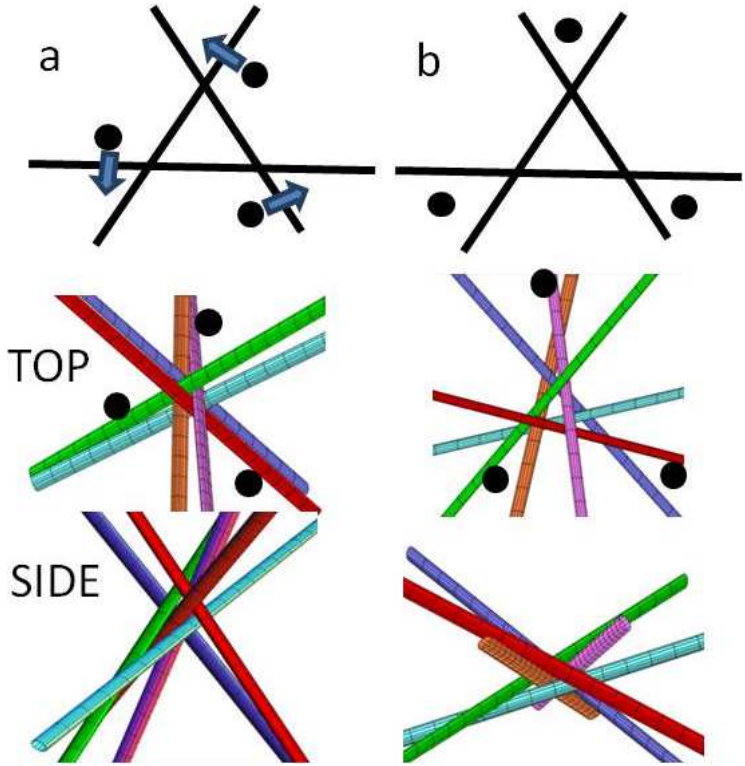}
\caption{6-knots with the chiral matrix of the determinant $-125$ and invariant $9.6666\ldots$. (a) Schematic, top, and side view of a unique configuration of the 6-knot present in quadruple in the 8-cross with the determinant 1625. (b) A slightly different 6-knot that cannot be present in the 8-cross. The difference is illustrated by the black disks and arrows indicating the switching between the cylinders.}
\label{fig:3}       % Give a unique label
\end{figure}

In Fig. \ref{fig:3}a,b we show illustrations of two $C_3$ symmetric 6-knots with equal elliptic cylinders (elliptic cylinders will be discussed in detail below) both with invariant $9.6666\ldots$ . These 6-knots are made of two sandwich triangles (called wreathes in \cite{Ref3}) rotated in opposite directions around $C_3$  axis. The "upper" wreath in a TOP view in Fig.\ref{fig:3}a,b rotates clockwise and consists of unknotted cylinders (1-red, 3-green, 5-violet) and "lower" wreath rotates counterclockwise and consists of knotted cylinders (0-brown, 2-blue, 4-cyan) which exactly corresponds to the Ring matrix given in Eq. (\ref{eq:12}).

The difference between 6-knots in Fig. \ref{fig:3}a and Fig. \ref{fig:3}b is emphasized by the black disk which marks the cylinders of the "upper" wreath: one can notice that pairs of cylinders (5th and 0th), (3rd and 4th), and (1st and 2nd) switch their chirality without changing the Ring matrix and the determinant $-125$ of the chirality matrix. The upper panels in Fig. \ref{fig:3} illustrate these two configurations more schematically which helps us below to analyze the case of 8-cross with the determinant $1625$.

The most important for the analysis with the schematic is the position of the black disks either in the sectors adjacent to the sides of the triangle for one configuration in Fig. \ref{fig:3}a or in the sectors adjacent to the vertexes of the triangle. Finally, we found that the configuration in Fig. \ref{fig:3}b switches into the one with the invariant 11.6842… (see Appendix 3) if one lifts the "upper" wreath father from the "lower" one, then the chirality matrix remains the same but new rings will be created passing a degenerate coplanar configuration. This explains why the configuration in Fig. \ref{fig:3}b is 1.33 times less frequent than the one in Fig. \ref{fig:3}a as it is clear from the results of the statistics of random configurations given in Appendix 2 where invariants for both configurations are highlighted in bold.

\begin{figure}[h]
  \includegraphics [width=0.85\textwidth]{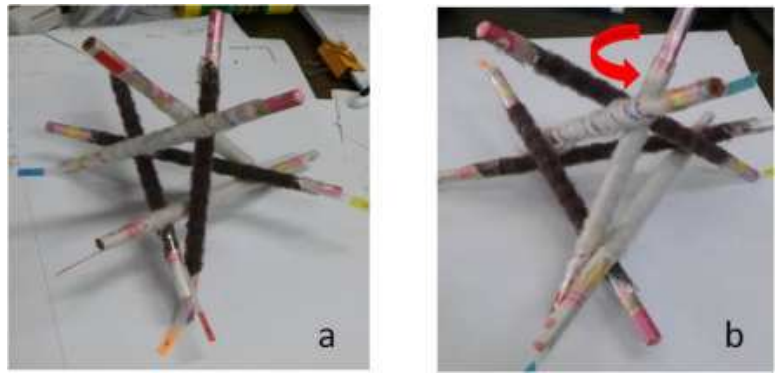}
\caption{A pencil model for 6-cross of (a) 9.666… with three free cylinders and (b) 5.8928571429 with all knotted cylinders.  The red arrow shows the direction of the rotation of the "upper" three cylinders in (a) to turn this configuration into (b).}
\label{fig:4}       % Give a unique label
\end{figure}

In order to demonstrate the amount of control gained over all 6-crosses by the invariants found let us show another 6-cross with $C_3$ symmetry and with a rather rare invariant of 5.8928571429 from Appendix 2 and 3. This 6-cross has all its cylinders knotted which its Ring matrix shows:
   \begin{equation}
R=\left( \begin{array}{cccccc}
0&1&3&1&1&3\\
3&0&1&3&1&1\\
1&3&0&1&3&1\\
1&1&3&0&1&3\\
3&1&1&3&0&1\\
1&3&1&1&3&0 \end{array} \right)
\label{eq:13}
\end{equation}
If one compares Eq. (\ref{eq:12}) with Eq. (\ref{eq:13}) one can see that there are no free cylinders in Eq. (\ref{eq:13}). The 6-cross 5.8928571429  can be easily produced from the one in Fig. \ref{fig:3}a (also given in Fig. \ref{fig:4}a as a pencil model) by rotating its "upper" sandwich triangle on 120 degrees. The resulting beautifully symmetric configuration, shown in Fig. \ref{fig:4}b has the chiral tetrahedral symmetry group. It is also interesting to note that the configuration can be decomposed into two "free" pieces only by moving a 3-cross away: one cylinder or two cylinders cannot be moved away being encaged. In this way it coincides with the famous Six Piece Wooden Star Puzzle.

In Appendix 3 there is the rarest 6-cross with the determinant $-125$ (it could not even be found with a statistic procedure) with invariant $5.2175438596$ which has the cylinders even more entangled than the one in Fig. \ref{fig:4}b, but we leave this subject for further publications.

\begin{figure}[h]
  \includegraphics [width=0.85\textwidth]{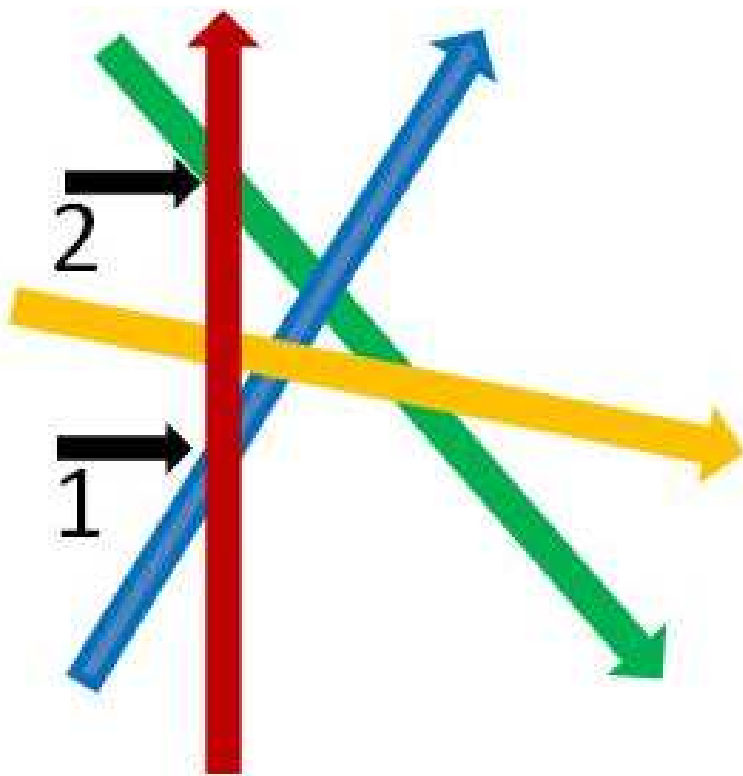}
\caption{A would-be 4-knot with two rings that illustrates the impossibility of mutually touching. The arrows mark two touch points that cannot be realized simultaneously.}
\label{fig:5}       % Give a unique label
\end{figure}

We made sure that all 6-crosses with the determinant $-125$ except the ones with the invariant $9.6666\ldots$ cannot be 6-knots.  Indeed, it turned out that all of them contain a special 4-cross that cannot be a 4-knot, that is the configuration cannot have all 4 cylinders in mutually touching. This special 4-cross is simply characterized by a Ring matrix that has two entangled cylinders and two free ones. Commonly, 4-crosses have either no entangled cylinders or only one entangled by the triangle wreath of the rest of the cylinders.  Having two entangled cylinders produces the highly restricted configuration which can be illustrated in Fig. \ref{fig:5}.

Indeed, in Fig. \ref{fig:5} the blue cylinder is inside the ring made by the red, yellow and green cylinders while the yellow cylinder is in the ring made by the red, blue and green cylinders, provided that all cylinders touch mutually. Yet, the simultaneous touching, say, at points 1 and 2 is impossible: either 1 or 2 may be a touching point, never both. This is a see-saw case similar to the one studied in \cite{Ref3} with respect to a 7-cross. Thus, any 4-cross sub-configuration if having 2 rings in their Ring matrix prevents any $n$-cross (that contains this 4-cross as a sub-configuration) from turning into an $n$-knot. Using this criterion, we inspected all the 6-cross configurations of Appendices 2,3 and found out that all of them except the ones with the invariant $9.6666\ldots$ are forbidden for mutually touching.

Finally, it is interesting to note that using the Ring matrix for selection procedure of realizable configurations produces other criteria. One to mention is a 5-cross with all free cylinders that is with zero Ring matrix. One can prove that such a 5-cross cannot be a 5-knot, because mutually touching demands at least one ring. Therefore if an $n$-cross configuration contains a 5-cross with zero Ring matrix then it cannot be an $n$-knot ever.
\section{\textbf{An 8-cross with determinant $1625$ cannot be an 8-knot} }
\label{sec:4}
Now it is left to prove that even if an 8-cross configuration with the determinant $1625$ contains four allowed 6-knots with the determinant $-125$ and invariant $9.6666\ldots$ still the 8-knot may be forbidden.
We constructed the 8-cross by using $C_4$ symmetry discussed above. To build-in the four 6-knots two schemes given below are possible. Some lines may look parallel but it is only due to projection. Black disks as in Fig. \ref{fig:3} mark the positions where the "upper" wreath of unknotted 4 cylinders punctures the "plane" of the "lower" wreath of opposite orientation. One can see that only scheme (a) can produce four 6-knots of $9.6666\ldots$ invariant. Indeed, if one removes a black disk together with the adjacent side line of the square one has exactly the schematic of Fig. \ref{fig:3}a - a triangle with black disks each within the sector of the side of the triangle. On contrary, if one removes a black disk with an adjacent line from scheme (b) then only two disks remain in the vertex sectors and one disk inevitably gets into the side sector which is different from Fig. 3b. Therefore such a configuration cannot produce four 6-knots needed.
\begin{figure}[h]
  \includegraphics [width=0.85\textwidth]{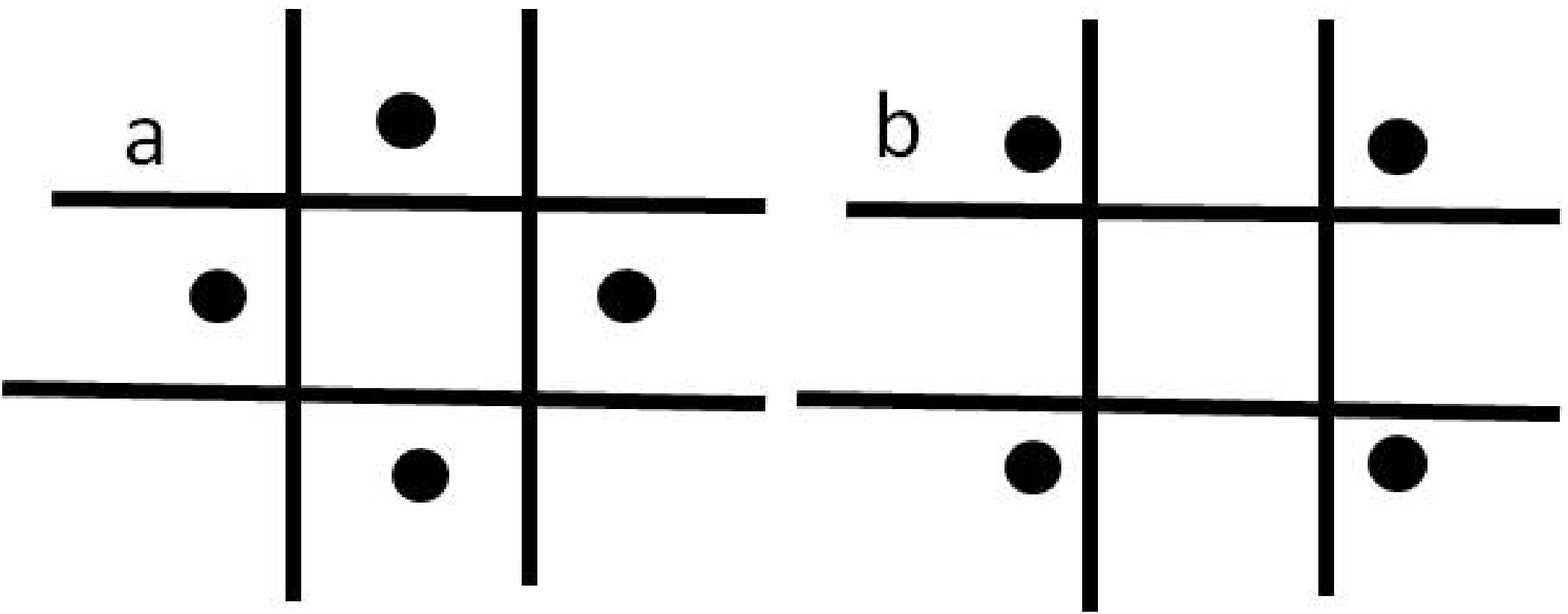}
%\caption{}
\label{fig:6}       % Give a unique label
\end{figure}
The configuration built on scheme (a) is characterized by the chirality matrix of Eq. (\ref{eq:9}) and the Ring matrix:
\begin{equation}
R=\left( \begin{array}{cccccccc}
0&0&0&0&0&0&0&0\\
5&0&7&7&3&5&3&3\\
0&0&0&0&0&0&0&0\\
3&3&5&0&7&7&3&5\\
0&0&0&0&0&0&0&0\\
3&5&3&3&5&0&7&7\\
0&0&0&0&0&0&0&0\\
7&7&3&5&3&3&5&0 \end{array} \right)
\label{eq:14}
\end{equation}

\begin{figure}[h]
  \includegraphics [width=0.85\textwidth]{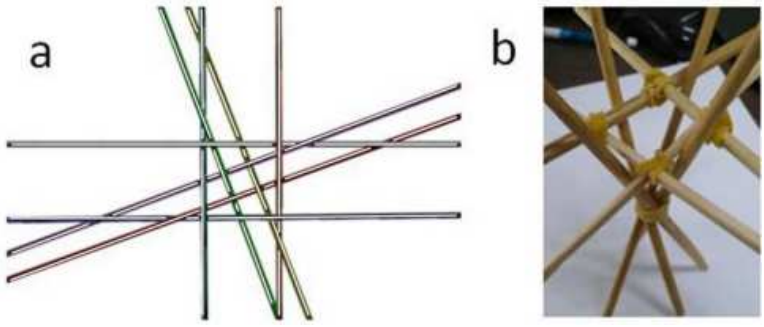}
\caption{(a) The top view of the 8-cross with the determinant $1625$; (b) its wooden model.}
\label{fig:7}       % Give a unique label
\end{figure}

The determinant of this chirality matrix is indeed $1625$ that can be considered as an "experimental" verification of the construction based on $C_4$-symmetry. The corresponding 8-cross is shown in Fig. \ref{fig:7} in two projections one of which in Fig. \ref{fig:7}a is the view along $C_4$ axis and the other is a general view of a chopstick model. Its invariants are $23.304029304$ and $25.7509157509$ for the mirror configuration.

This $P1625$ 8-cross with four sub-configuration 6-knots with the invariant $9.6666\ldots$ and the determinant $-125$  should be unique (provided that we indeed exhaust and control numerically all configurations of the determinant $-125$ discussed above): if one reduces the symmetry to, say, $C_2$ then four 6-crosses with the determinants $-125$ and with two pairs of different invariants could only exist as sub-configurations, for instance, two $9.6666\ldots$  and two $10.52$ etc. When four $9.6666\ldots$ reappear accidentally, then the symmetry is restored to $C_4$ and the 8-cross is the unique configuration.

Yet, the inspection of the unique 8-cross revealed the presence of a 4-cross with two rings as a sub-configuration. This makes the total 8-cross configuration unable to be turned into an 8-knot, as far as a 4-cross with two rings can never be turned into a 4-knot as we proved before. Thus we proved that all possible 8-crosses with the chirality matrix $P1625$ are forbidden as 8-knots.
Note that generally the symmetry here works against $n$-knots because symmetry imposes more restrictions and tacitly reduces the number of the degrees of freedom. This situation is paradoxical and different from common situations in geometry when symmetry favors solutions.

It immediately follows that an 11-cross with $M11$ as a chirality matrix is impossible as an 11-knot for any cylinders of arbitrary cross-section because the configuration should contain $P1625$. It reduces the number of cylinders in $n$-knots down to 10 instead of previously established upper limit of 14 cylinders \cite{Ref3}. This is the main result of our paper. Below we show that 10-knots are realizable at least with unequal elliptic cylinders.

\section{	\textbf{$n$-knots with elliptic cylinders}}
\label{sec:5}

After elucidating topological matters let us turn to the metric questions. First of all, to be able to build a 10-knot one cannot use round cylinders anymore, even those of arbitrary radii because the degrees of freedom for round cylinders are exhausted: there are too many restrictions of mutual contacts. The maximum number of round cylinders for an n-knot is $n=9$ as it was proven in \cite{Ref2}. Yet, it is enough to increase the number of the degrees of freedom of individual cylinders by using the cylinders with elliptical cross-section. It opens many interesting possibilities. For example, in \cite{Ref3} it was shown that there is only one 7*-knot (plus its mirror image) that has all equal round cylinders. Numerous other 7-knots that we found can exist only with cylinders of unequal radii. Meantime, 7-knots with all equal elliptical cylinders which were not possible to make of all equal round cylinders exist.  Define such configuration as 7**-knot, probably all 7-knots with unequal round cylinders can be turned into 7**-knots.

Otherwise, any elliptical 7-knot with arbitrary elliptical cylinders can be "rounded", that is deformed into a 7-knot with round cylinders. The same situation holds for 8-knots: we succeeded in building 8*-knot configurations with all equal elliptic cylinders and in "rounding" them as well. However, a 9-knot is different: although it is still possible to make a 9*-knot with all equal elliptic cylinders (actually it looks more like the one with equal thin stripes), there are elliptic configurations that cannot be "rounded". That is some elliptic configurations do not exist for round cylinders at all. We will discuss the difference below. As to 10-knots we have not found a configuration with all equal elliptic cylinders though it would be very interesting to find it because such a configuration would exhaust all the degrees of freedom as we will see below.
The system of equations for elliptic cylinders introduces more variables compared to round cylinders. Indeed, for round cylinders we have 5 variables for each cylinder \cite{Ref5}:

$t_i,p_i$ -  two angles in spherical coordinates that define the unit vector of the direction of the $i$th cylinder axis $\vec n_i=\vec n(t_i,p_i)$ ; $x_i,y_i$ -  two coordinates of the point where the axis punctures the plane $xy$ that define the plane vector $\textbf{v}_i=(x_i,y_i,0)$; finally, the radius $r_i$  of the  cylinder.

Note that one can turn an $n$-cross into its mirror configuration by just replacing all  $\textbf{v}_i$ into $-\textbf{v}_i$ while preserving all orientations. Of course such an operation needs multiple crossing the lines. The operation is analogous to turning a glove inside-out to transform a left glove into the right one.

For elliptic cylinders we have 7 variables: instead of  $r_i$ there are two semi-axes of the ellipse $a_i,b_i$ and the Euler angle $\omega_i$  of its rotation around the cylinder axis. This is more than enough for constructing a 10-knot. However, if we wanted to find a configuration with all equal elliptic cylinders then we would find again 5 variables for each cylinder (instead of  $r_i$ for round cylinders it would be  $\omega_i$   and the  aspect ratio $b/a$ of semi-axes) . $a$ itself cannot be considered as an independent variable because it can be used for scaling the distances. Then the total degrees of freedom \cite{Ref5} for $n=10$ all equal elliptic cylinders would be completely exhausted by the pairwise mutual contacts and 6 degrees of the solid body:
\begin{equation}
5n+1-n(n-1)/2-6=0
\label{eq:15}
\end{equation}

\begin{figure}[h]
  \includegraphics [width=0.85\textwidth]{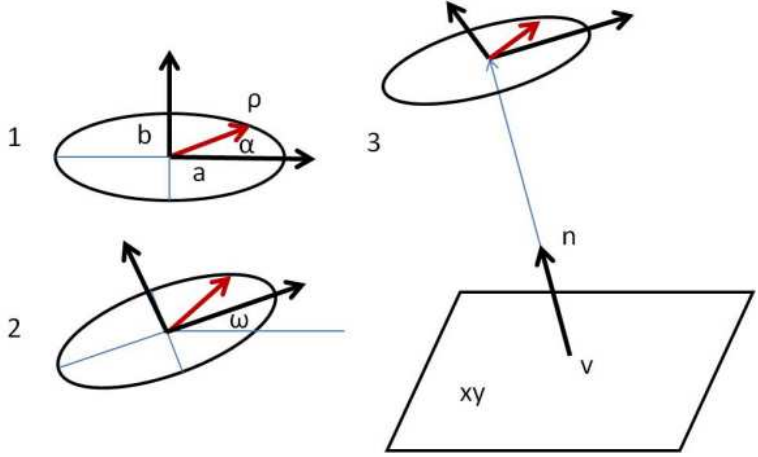}
\caption{Schematic of the stepwise construction of an elliptical cylinder cross-section position.}
\label{fig:8}       % Give a unique label
\end{figure}

The geometry of the elliptic cylinder is illustrated in Fig. \ref{fig:8}. To explain it in the simplest way let us construct a general position of the elliptic cylinder by steps.

	 1.Take the ellipse which major semi-axis of the size $a$ lies along $x$ and minor semi-axis of the size $b$ lies along $y$. It can be obtained by the affine deformations of the unit circle along the axes. The radius-vector ${\vec\rho}$ to a point on the ellipse allows parametrization by the angle $\alpha$ between the radius vector of the unit circle and $x$-axis unit vector $\textbf{n}_a=(1,0,0)$ which is
${\vec\rho}=a\textbf{n}_a  \cos \alpha + b\textbf{n}_b  \sin \alpha$
(vector $\textbf{n}_b=(0,1,0)$).

	 2. Rotate the ellipse around $z$ on the angle  $\omega$ so that the unit vector along the major semi-axis becomes $\textbf{n}_a=(\cos\omega,\sin\omega,0)$ and the one along the minor semi-axis becomes $\textbf{n}_b=(-\sin\omega,\cos\omega,0)$ .

	 3. Rotate the ellipse with the rotation matrix defined by the spherical angles:

$$
Rot(p,t)=\left( \begin{array}{ccc}
\cos p&-\sin p&0\\
\sin p&\cos p&0\\
0&0&1 \end{array} \right) \left( \begin{array}{ccc}
\cos t&0&\sin t\\
0&1&0\\
-\sin t&0&\cos t\end{array} \right)
\label{eq:15a}
$$

so that
\begin{equation}
 \vec \rho=a\textbf{N}_a  \cos\alpha+ b\textbf{N}_b  \sin\alpha
\label{eq:16}
\end{equation}

where $\textbf{N}_a= Rot(p,t) \textbf{n}_a$ and $\textbf{N}_b= Rot(p,t) \textbf{n}_b$.

Note that our consideration is not restricted with the specific form of Eq. (\ref{eq:16}). In principle, ${\vec \rho}$  can be any smooth function on $\alpha_i$ to draw the cross-section of the $i$th cylinder where we equipped $\alpha$ with the subscript $i$. If smoothness is violated at some isolated points then broader "non-holonomic" possibilities for contacts appear, up to the infinite number of stripes (as degenerate ellipses) touching at one point.

Now we have to bring into a contact two elliptic cylinders. At the point of contact, tangent vectors $\partial \vec\rho_i/\partial \alpha_i $ for ellipses of both cylinders should lie in the plane orthogonal to vector  $[\textbf{n}_i\times \textbf{n}_j ]$  thus satisfying conditions:

\begin{equation}
\partial {\vec\rho}_i/\partial \alpha_i  [\textbf{n}_i\times \textbf{n}_j ]=0
\label{eq:17}
\end{equation}

If one uses Eq. (\ref{eq:16}) then one obtains from Eq. (\ref{eq:17}) $\tan \alpha_i=(b/a)  (\textbf{N}_b  [\textbf{n}_i\times \textbf{n}_j ])/(\textbf{N}_a  [\textbf{n}_i\times \textbf{n}_j ] )$ and, correspondingly, finds ${\vec\rho}_i$ from Eq. (\ref{eq:16}).
 When viewed along the plane orthogonal to bivector $[\textbf{n}_i\times \textbf{n}_j ]$, the cylinders look like two parallel stripes touching one another along the edge line. That means that the half-widths of the projected stripes sum up to the distance between the projected axes of the cylinders. It is easy to understand that this distance between the axes on the projection is  $|[\textbf{n}_i\times \textbf{n}_j ](\textbf{v}_i-\textbf{v}_j ) |$.
Let us write down the contact equations for two touching elliptic (or another shape for arbitrary  ${{\vec\rho}}_i(\alpha_i)$) cylinders:

\begin{equation}
|[\textbf{n}_i\times \textbf{n}_j ]{\vec\rho}_i|+|[\textbf{n}_i\times \textbf{n}_j ]{\vec\rho}_j |=|[\textbf{n}_i\times \textbf{n}_j ](\textbf{v}_i-\textbf{v}_j ) |
\label{eq:18}
\end{equation}

The equation for round cylinders would simplify:
\begin{equation}
      (r_i+r_j )|[\textbf{n}_i\times \textbf{n}_j ]|=|[\textbf{n}_i\times \textbf{n}_j ](\textbf{v}_i-\textbf{v}_j ) |
 \label{eq:19}
 \end{equation}

Eqs. (\ref{eq:18}) and (\ref{eq:19}) could be re-written without absolute values on the rhs with the help of the chirality matrix:

\begin{equation}
{ P}_{i,j} (|[\textbf{n}_i\times \textbf{n}_j ]{\vec\rho}_i|+|[\textbf{n}_i\times \textbf{n}_j ]{\vec\rho}_j | )=[\textbf{n}_i\times \textbf{n}_j ](\textbf{v}_i-\textbf{v}_j )
 \label{eq:20}
 \end{equation},
\begin{equation}
{ P}_{i,j} (r_i+r_j )|[\textbf{n}_i\times \textbf{n}_j ] |=[\textbf{n}_i\times \textbf{n}_j ](\textbf{v}_i-\textbf{v}_j )
 \label{eq:21}
 \end{equation}

On the rhs we have the non-normalized chirality matrix entries  ${\cal P}_{i,j}$ discussed above (see Eq. (\ref{eq:3})). This induces the strong restrictions to the determinants of the touching cylinders to be zero when the rank is above 6.
Eq. (\ref{eq:20}) may give a broader number of solutions than Eq. (\ref{eq:21}) for round cylinders. That explains why we found a solution for 9-knots that does not have its counterpart with round cylinders.

As we mentioned before, for $n>6$ the determinant of the non-normalized chirality matrix Eq. (\ref{eq:3}) which coincides with the rhs of Eqs. (\ref{eq:20}), (\ref{eq:21}), is zero.  This fact is especially meaningful for $n=9$ and Eq. (\ref{eq:21}) for round cylinders. Indeed, the number of zero determinant sub-matrices for $n=9$ is in total 46: one comes from the determinant of the matrix itself, nine of them come from zero determinants for $(8\times 8)$ sub-matrices and 36 come from $(7\times 7)$ sub matrices. This number exceeds the number of degrees of freedom of a 9-knot. Thus some strong selection rule may be imposed that selects out chirality matrices in Eq. (\ref{eq:21}). So many zero determinants may explain the general tendency in 9-knot chirality matrices to have small determinants; most frequently for 9-knots that we found, they are of zero determinant.

The situation gets even more dramatic for a 10-knot: the number of zero equations becomes 386. We leave these questions for the future work.  Here we only mention that we arranged zero determinants for matrices of rank 9 into a Mathcad11 program that tries to satisfy all the conditions starting from random initial values for angles and radii and with a given chirality matrix.  We tried how it works as a selection rule on the matrix of rank 9 extracted from $M11$ (see Eq. (\ref{eq:8})) which contains the matrix $P1625$ discussed above. Indeed, no solution was found by the Mathcad11 program for this matrix though it easily gave solutions for non-forbidden matrices of rank 9.

Eqs. (\ref{eq:20}), (\ref{eq:21}) could be interpreted as follows: the "product" of the spatial plane and the "Riemann plane" on the rhs of Eqs. (\ref{eq:20}), (\ref{eq:21}) is kept "constant" (as the dimensions or radii on the lhs are constants). Thus, Eqs. (\ref{eq:20}),(\ref{eq:21}) relate a pure "geometry matrix" that contains space vectors on the rhs, with the "matter matrix" on the lhs, resembling field equations with the source. Indeed, the source bunches free lines which from a large distance look like radially divergent field lines issuing from the source. The reverse problem of defining the structure of the source while being given the lines far from the source can be posed as well.

\begin{figure}[h]
  \includegraphics [width=0.85\textwidth]{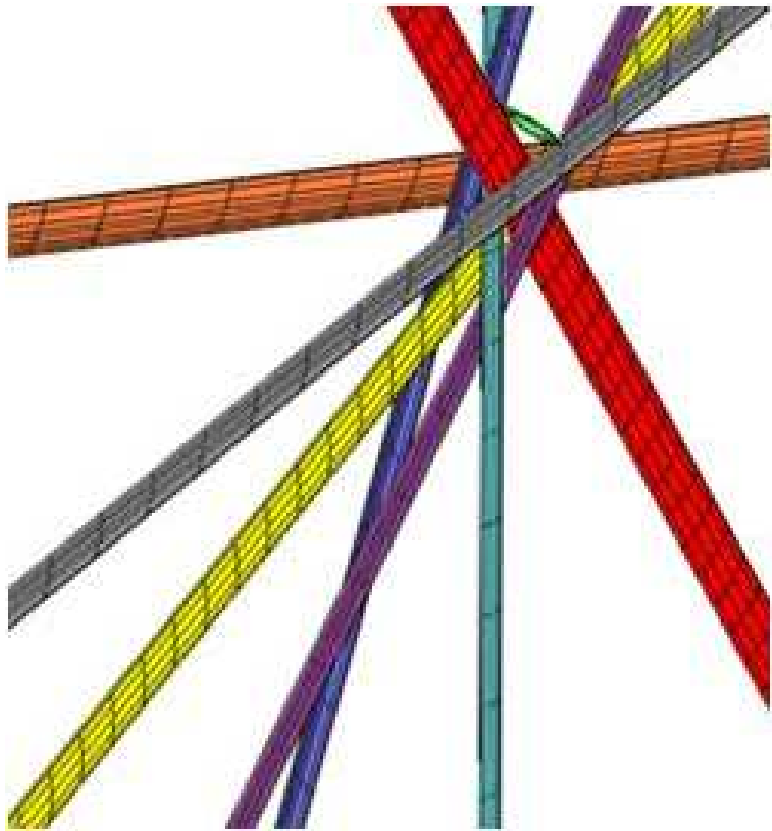}
\caption{An 8*-knot of all equal elliptic cylinders. We deliberately projected the 8*-knot so that one can see the elliptic cross-section of one of the cylinders.}
\label{fig:9}       % Give a unique label
\end{figure}

\begin{figure}[h]
  \includegraphics [width=0.85\textwidth]{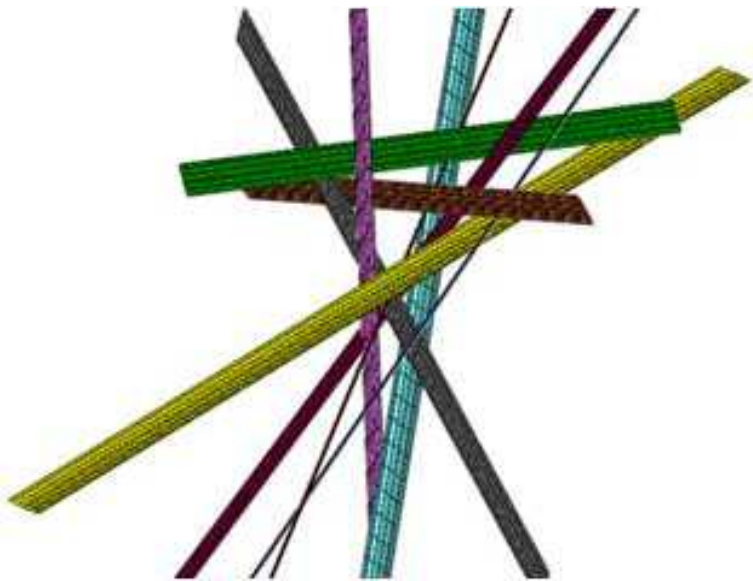}
\caption{An 9*-knot of all equal elliptic cylinders. The aspect ratio is rather large therefore the equal elliptic cylinders look more like equal stripes (certainly, the projection makes them look unequal in the figure).}
\label{fig:10}       % Give a unique label
\end{figure}

We solved Eq. (\ref{eq:18}) for 8, 9 and 10-knots with Mathcad11 in the same way as we described solving Eq. (\ref{eq:19}) in our previous works. Here we give a brief account of the results showing an 8*-knot with all equal elliptic cylinder of aspect ratio 0.15 (Fig. \ref{fig:9}) with the invariant 20.4666666667, mirror invariant 20.9619047619, and the determinant of the chirality matrix -7,
and a 9*-knot (Fig. \ref{fig:10}) with all equal elliptical cylinders with the aspect ratio of 0.025. The invariant of the latter is 15.9676775014 and its mirror invariant is 40.2095514799 calculated according to Eq. (\ref{eq:6}) with the chirality matrix and the Ring matrix given below:
     \begin{equation}
P=\left( \begin{array}{rrrrrrrrr}
0&-1&-1&-1&-1&-1&-1&-1&-1\\
-1&0&-1&-1&1&-1&-1&1&-1\\
-1&-1&0&1&1&-1&-1&-1&1\\
-1&-1&1&0&-1&-1&1&-1&1\\
-1&1&1&-1&0&-1&1&1&1\\
-1&-1&-1&-1&-1&0&1&-1&1\\
-1&-1&-1&1&1&1&0&-1&1\\
-1&1&-1&-1&1&-1&-1&0 &1\\
-1&-1&1&1&1&1&1&1 &0 \end{array} \right)
\label{eq:22}
\end{equation}

    \begin{equation}
R=\left( \begin{array}{ccccccccc}
0&2&2&10&2&5&5&2&2\\
8&0&7&8&7&8&7&8&7\\
0&0&0&0&0&0&0&0&0\\
0&0&0&0&0&0&0&0&0\\
0&0&0&0&0&0&0&0&0\\
0&0&0&0&0&0&0&0&0\\
3&3&6&3&9&9&0&3&6\\
4&3&12&4&3&3&3&0 &4\\
6&6&6&10&10&5&5&6&0 \end{array} \right)
\label{eq:23}
\end{equation}

The determinant of the chirality matrix of Eq. (\ref{eq:22}) is zero.

Most of all sub-configurations of seven mutually touching cylinders that can be extracted from 8*-knots or 9*-knots are 7**-knots with all equal elliptical cylinders which were not possible to make of all equal round cylinders. All 8-knots as well as many 9-knots can have their round counterparts with the same topological characteristics described above and given in \cite{Ref3}. However, we found a specific 9-knot that does not have its round counterpart which is in agreement with our earlier statement that Eq. (\ref{eq:18}) has larger set of solutions compared to Eq. (\ref{eq:19}).
\begin{figure}[h]
  \includegraphics [width=0.85\textwidth]{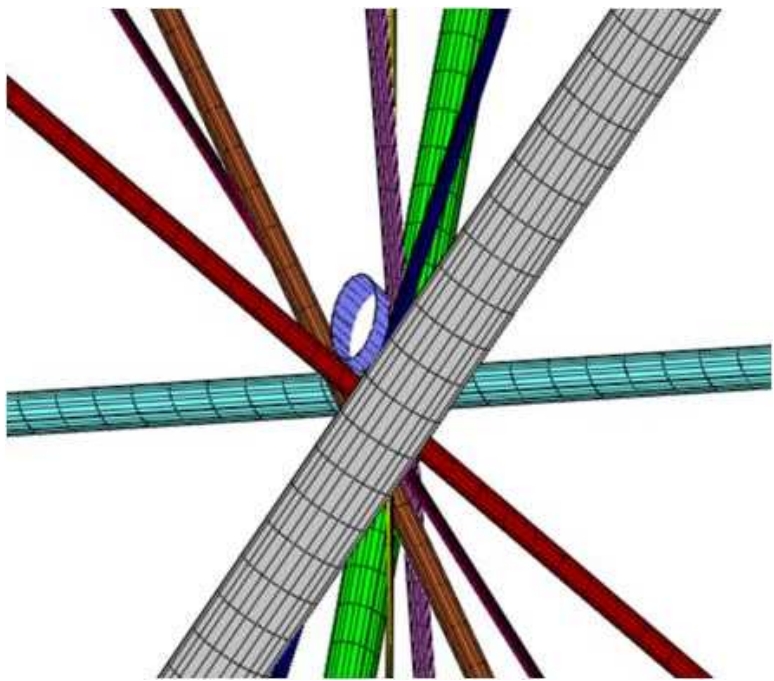}
\caption{A 10-knot of elliptic cylinders.}
\label{fig:11}       % Give a unique label
\end{figure}

As for 10-knots which cannot have solutions with round cylinders, we found several configurations by solving Eq. (\ref{eq:18}); one is presented in Fig. \ref{fig:11}.
Its invariants are $34.7078032583$ and $17.7482857444$ for the mirror one.  The determinant of the chirality matrix is $-25$.

This 10-knot has the parameters listed in the tables below:
\begin{table}[h]
\caption{Computationally obtained parameters for a 10-knot.}
\label{tab:1}       % Give a unique label
\begin{tabular}{lllll}
\hline\noalign{\smallskip}
$t_i$ & $p_i$	& $x_i$	& $y_i$ & $\omega_i $ \\
\noalign{\smallskip}\hline\noalign{\smallskip}
0	& 0	& 0	& 0	& 3.1254845844 \\
0.8888266921& 	4.8105113254& 	1.718837714	& -1.0154516939	& 3.0137625571\\
1.7307454079& 	4.3698637441& 	15.8387596901& 	35.4191212761& 	3.798746922\\
0.6710427219& 	2.1610133993& 	0.4876982243& 	-4.7627389073& 	2.9915865287\\
1.9722582511& 	3.3409580456& 	6.3213598491& 	-0.9629547593& 	2.460123252\\
0.5053571239& 	3.5602584953& 	3.410770679	& 0.0642137547	& 3.1415926536\\
2.0636231502& 	-0.0462841434& 	4.9265275262& 	3.7039199491& 	3.1391795035\\
2.7052117458& 	-0.3157789613& 	2.3817940085& 	0.3228529593& 	3.1415926536\\
2.3433579644& 	4.2721078026& 	0.3516109986& 	-2.2031892116& 	1.4094265204\\
1.1462686715& 	3.7874981154& 	-1.7994972106& 	-3.0922913208& 	1.9871990245\\
\noalign{\smallskip}\hline
\end{tabular}
\end{table}

\begin{table}[h]
\caption{Computationally obtained cylinder semi-axes for a 10-knot.}
\label{tab:2}       % Give a unique label
\begin{tabular}{ll}
\hline\noalign{\smallskip}
$a_i$& $b_i$ \\
0.8756395562& 	0.8756395562\\
1.1107047384& 	0.7276292249\\
2.989157107& 	1.4262671635\\
1.8316025399& 	1.3574498777\\
1.8438264359& 	0.8079920149\\
0.8219266058& 	0.4523138356\\
3.0522525128& 	3.0522524955\\
0.2147737114& 	0.1709326871\\
0.3864418576& 	0.010000001\\
0.5582558938& 	0.0100000003\\
\noalign{\smallskip}\hline
\end{tabular}
\end{table}
The chirality and Ring matrices are:
  \begin{equation}
P=\left( \begin{array}{rrrrrrrrrr}
0&1&1&1&1&1&-1&1&1&-1\\
1&0&1&1&1&-1&-1&1&-1&-1\\
1&1&0&-1&-1&-1&-1&1&-1&1\\
1&1&-1&0&-1&1&-1&-1&-1&1\\
1&-1&-1&-1&0&1&1&-1&1&-1\\
1&-1&-1&1&1&0&-1&-1&-1&-1\\
-1&-1&-1&-1&1&-1&0&1&1&1\\
1&1&1&-1&-1&-1&1&0 &-1&-1\\
1&-1&-1&-1&1&-1&1&-1 &0&-1\\
-1&-1&1&1&-1&-1&1&-1&-1 &0 \end{array} \right)
\label{eq:24}
\end{equation}
and
    \begin{equation}
R=\left( \begin{array}{cccccccccc}
0&4&4&10&6&4&14&8&4&6\\
7&0&11&7&9&11&9&11&7&9\\
0&0&0&0&0&0&0&0&0&0\\
0&0&0&0&0&0&0&0&0&0\\
2&2&12&6&0&2&2&2&6&2\\
0&0&0&0&0&0&0&0&0&0\\
0&0&0&0&0&0&0&0&0&0\\
4&4&6&4&4&14&10&0 &8&6\\
9&5&5&15&5&9&5&5&0&5\\
8&6&14&6&12&6&8&6&6&0 \end{array} \right)
\label{eq:25}
\end{equation}

\section{Conclusion}
\label{sec:6}
By constructing 10-knots we reached the limit for the number of mutually touching infinite straight cylinders which is 10 as we proved. While solving the question we developed several tools to get more control on topology and geometry of configurations of straight infinite cylinders of arbitrary cross-section in 3D. Some of the tools need much more elaboration to become useful for any other applications. Generalizations may be foreseen for the chirality and Ring matrices (as well as for the invariant) to simplify the description of $n$-crosses and $n$-knots, their combinations and their sub-structures.

%\begin{acknowledgements}

%\end{acknowledgements}



\begin{appendix}
\section{\textbf{Appendix 1}}

The proof that any chirality (Seidel's adjacency) matrix of rank more than 18 contains $K5$ as a diagonal sub-matrix.

Let us take a $P19$ --- a chirality matrix of rank 19. Transform $P19$ to $T P19 T$ with the help of matrix $T$ that has all zero entries except the diagonal where $T_{00}=1$ and $T_{ii}=P19_{0i}$, where $i$ takes the values from 1 to 18. Then the transformed matrix will have all 1s in zero row (and column) except on the diagonal.

A submatrix of $P19$ obtained by elimination of the zero row and zero column (we enumerate the rows and columns starting from number zero), is of rank 18. Denote it as $P18$. $P18$ shall contain a diagonal submatrix $K4$ (its entries are either all +1's or -1's except all zeroes on the diagonal). Indeed, from the well-known fact that Ramsey's number $R(4,4)=18$ \cite{Ref9} it follows that the graph that corresponds to Seidel's matrix $P18$ (or its complementary graph) shall contain a complete graph of 4 nodes, therefore its adjacency matrix must contain the abovesaid diagonal submatrix $K4$. If non-zero entries of $K4$ are +1 then along with +1's of zero row and column of $P19$ it makes $K5$. If the non-zero entries of $K4$ are -1, multiply zero row and column of $P19$ by -1 which action does not change $K4$ from $P18$. By elimination all rows and columns from $P19$ except number zero ones and those of $K4$ we obtain a diagonal submatrix $K5$ which entries are all +1's or -1 except zeroes on the diagonal.

\section{\textbf{Appendix 2}}
First column gives the invariant for a configuration of 6-cross with the determinant of the chirality matrix -125 calculated according to Eq. (6). Out of trial number 49233 of configurations of mutually contacting cylinders, obtained from the set of randomly oriented axes of cylinders. the second column shows the frequency of a given configuration. The third column is 0 when there is no 6-knot and 1 otherwise.  The most probable and the only one realizable as a 6-knot is given in bold.
The fourth column shows the new invariant that discriminates all the configurations and their mirror ones (the fifth column). The sixth column gives the sum of the fourth and the fifth columns to demonstrate the coincidence when there is the same Ring matrix as for 9.6666… invariant.
\setcounter{table}{0}
\renewcommand{\thetable}{A\arabic{table}}

%\begin{table}[h]
\begin{longtable}{llllll}
%\caption{}
\label{tab:1}       % Give a unique label
%\begin{tabular}{llllll}
%\hline\noalign{\smallskip}
\\
8.33333	 &            3472&	0&	-0.86349&	-0.85898&	-1.72247\\
10.52	&             3253&	0&	-0.35797&	-0.88583&	-1.2438\\
8.86667	&             3117&	0&	-0.65779&	-0.71131&	-1.36909\\
8.82857	&             2499&	0&	-0.63005&	-0.83023&	-1.46027\\
\textbf{9.66667}&	             \textbf{2222}&	\textbf{1}&	\textbf{-0.32718}&	\textbf{-0.73902}&	 \textbf{-1.06621}\\
8.88235	&             2218	&0&	-0.64839&	-0.5668	&-1.21519\\
10.33333&	             1770&	0&	-0.53772&	-1.68779&	-2.22551\\
11.17391&	             1687&	0&	-1.04194&	-0.39432&	-1.43626\\
\textbf{9.66667}	&             \textbf{1618}&	\textbf{1}&	\textbf{-0.54683}&	\textbf{-0.51938}&	 \textbf{-1.06621}\\
10.06977&	             1409&	0&	-1.22533&	-0.3517&	-1.57704\\
9.17949	&             1211&	0&	-1.5472&	-0.45858&	-2.00578\\
8.88235	&             1156&	0&	-0.8905	&-0.32002&	-1.21052\\
10.33333&	984&	0&	-2.18256&	0.40353&	-1.77903\\
7.33333&	947	&0&	-0.9933&	-1.1739	&-2.1672\\
10.66667&	924	&0&	-0.39193&	-0.90613&	-1.29806\\
9.01869&	920	&0&	-1.14795&	-0.34392&	-1.49187\\
9.23529	&904&	0&	-0.25557&	-2.06643&	-2.322\\
11.8&	878&	0&	-1.60728&	-0.28659&	-1.89388\\
8.89286&	876&	0&	-0.84713&	-0.43703&	-1.28416\\
10.40678&	874	&0	&-1.35318&	-0.47658&	-1.82976\\
8.41844&	778&	0&	-0.76749&	-1.29005&	-2.05753\\
9.78947&	656&	0&	-0.45231&	-0.95441&	-1.40672\\
9.8&	639&	0&	-0.54882&	-1.66044&	-2.20926\\
9.38&	588&	0	&-1.04248&	-0.58219&	-1.62467\\
7.51613&	586&	0&	-0.61054&	-0.67558&	-1.28612\\
8.84946&	586	&0	&-0.31222&	-1.16555&	-1.47777\\
7.09091	&580&	0&	-0.53887&	-0.77592&	-1.31479\\
10.22222&	580&	0&	-1.29533&	-0.9987&	-2.29403\\
11.89091&	575&	0&	-0.46372&	-1.34614&	-1.80986\\
10.05&	520	&0&	0.61015	&-3.06305&	-2.45289\\
10.93548&	487&	0&	-1.37187&	-0.37385&	-1.74572\\
10.8&	469&	0&	-0.62081&	-1.7619	&-2.38271\\
8.38272	&452&	0&	-0.835	&-0.78913&	-1.62413\\
10.02703&	451&	0&	-0.29835&	-0.97454&	-1.27289\\
8.60731&	449&	0&	-0.41685&	-1.34449&	-1.76134\\
7.78261&	441&	0&	0.2592&	-3.00002&	-2.74083\\
10.04&	425&	0&	-2.7237&	0.36133&	-2.36236\\
9.09091	&399&	0	&-0.74843&	-0.33833&	-1.08676\\
10.33333&	357	&0	&-0.37794&	-1.84757&	-2.22551\\
10.21053&	338	&0&	-1.15619&	-1.38863&	-2.54482\\
8.86139&	310&	0&	-0.41699&	-1.09011&	-1.5071\\
8.99024	&305&	0&	-1.15442&	-0.70831&	-1.86273\\
9.4	&303&	0&	-1.00218&	-0.56564&	-1.56782\\
8.87879&	263&	0&	-0.52924&	-1.78596&	-2.3152\\
11.21176&	248&	0&	-0.50986&	-1.45743&	-1.96729\\
11.68421&	239	&0&	-0.45924&	-1.15869&	-1.61793\\
9.42857	&215&	0&	-2.09315&	-0.59171&	-2.68486\\
7.22222&	206	&0	&-0.69849&	-0.46761&	-1.1661\\
8.52466&	195	&0	&-0.70257&	-1.08478&	-1.78735\\
8.41333&	178	&0	&0.0271	&-2.71597&	-2.68887\\
9.5	&167&	0	&-0.97545&	-0.73952&	-1.71497\\
9.1	&163	&0	&-0.3205&	-0.96304&	-1.28355\\
10.2&	156	&0&	-2.13564&	-0.75036&	-2.886\\
10.16667&	155&	0&	-0.32411&	-2.72199&	-3.0461\\
7.03023&	147&	0&	-0.85903&	-0.81604	&-1.67507\\
12.2&	145	&0	&-1.10625&	-0.52191&	-1.62816\\
8.73593	&141&	0&	-0.53393&	-1.57011&	-2.10404\\
9.24837	&140&	0	&-0.89254&	-0.53277&	-1.42531\\
9.97949&	139	&0&	-0.05771&	-2.09774&	-2.15545\\
11.93333&	137	&0&	-0.61361&	-1.30849&	-1.9221\\
9.2&	132&	0&	-0.53917&	-0.97056&	-1.50973\\
8.64602&	128&	0&	-1.12845&	-0.94358&	-2.07203\\
13.4&	121&	0&	21.81228&	11.28091&	33.09318\\
7.68116&	121	&0&	-0.72253&	-0.8494&	-1.57193\\
10.13333&	115	&0	&-0.51647&	-0.92326&	-1.43974\\
8.53271&	92&	0&	-0.70603&	-1.05362&	-1.75965\\
6.76119&	85	&0&	-0.61503&	-0.71285&	-1.32787\\
8.85185	&85&	0&	-0.62865&	-1.82022&	-2.44887\\
6.84426&	83&	0&	-0.80235&	-0.47962&	-1.28197\\
9.75088&	72&	0&	-1.14318&	-0.60945&	-1.75263\\
8.51055&	71&	0&	-0.7566	&-1.69047&	-2.44707\\
10.26667&	67	&0&	-0.81228&	-0.82276&	-1.63504\\
6.9&	54&	0&	-0.88872&	-0.53337&	-1.42209\\
8.85185&	54&	0&	-0.3126&	-2.13627&	-2.44887\\
8.93035&	48&	0&	-1.33215&	-0.40578&	-1.73793\\
6.92683&	46&	0&	-0.55653&	-0.6511&	-1.20763\\
8.28696&	42&	0&	-1.15792&	-0.40713&	-1.56505\\
10.03529&	40&	0&	0.23607&	-2.75072&	-2.51465\\
8.41176&	39&	0&	-1.06799&	-0.42837&	-1.49636\\
5.89286&	39&	0&	-0.63931&	-0.63931&	-1.27862\\
9.25352&	37&	0&	-1.11494&	-0.52475&	-1.6397\\
6.46552&	36&	0&	-0.63529&	-0.80835&	-1.44365\\
7.06818&	33&	0&	-0.44975&	-0.85016&	-1.2999\\
7.78571&	30&	0&	-0.47383&	-0.67468&	-1.14851\\
9.95094&	29&	0&	-1.92015&	-0.32179&	-2.24194\\
8.61905&	27&	0&	-1.4815&	-0.56487&	-2.04637\\
8.42424&	25&	0&	-0.86261&	-1.9542&	-2.81681\\
9.46897&	25&	0&	-1.28568&	-0.69988&	-1.98557\\
8.87571&	25&	0&	-1.29741&	-0.5374&	-1.83481\\
10.6&	23&	0&	-1.42679&	-1.00529&	-2.43208\\
8.66667&	22&	0&	-0.69507&	-1.46279&	-2.15785\\
9.04762&	17&	0&	-1.72501&	-0.83626&	-2.56127\\
8.10345&	17&	0&	-1.62937&	-1.57084&	-3.2002\\
7.54762&	15&	0&	-1.03887&	-0.7428&	-1.78167\\
7.40421&	14&	0&	-0.61253&	-1.11047&	-1.723\\
7.2619&	14&	0&	-0.64727&	-0.96147&	-1.60874\\
9.23242	&12&	0&	-1.66205&	-0.54407&	-2.20612\\
10.75&	11&	0&	16.14509&	27.48443&	43.62951\\
8.19048	&10&	0	&-1.18639&	-0.50581&	-1.6922\\
6.97619&	9&	0&	-0.72874&	-0.61815	&-1.34689\\
9.25&	7&	0&	-0.41403&	-1.13903&	-1.55306\\
7.80952&	7&	0&	-1.07492&	-0.40108&	-1.476\\
9.31206&	5&	0&	-1.34004&	-0.53945&	-1.87949\\
9.91111&	1&	0&	-1.62382&	-0.6921&	-2.31592\\
3.67925&	1&	0&	-0.49998&	-0.5281&	-1.02808\\
%\end{tabular}
%\end{table}
\end{longtable}

\section{\textbf{Appendix 3}}

Invariants of 6-knots with $C_3$ symmetry and the determinant of the chirality matrix -125.
The third and the fourth columns show the new invariants for the direct and mirror configurations.
\begin{table}[h]
\begin{tabular}{llllll}
7.33333&	-125&	-1.1739	&-0.9933\\
10.2&	-125&	-0.75036&	-2.13564\\
\textbf{9.66667}	&-\textbf{125}&	\textbf{-0.73902}&	\textbf{-0.32718}&\\
5.89286&	-125&	-0.63931&	-0.63931\\
\textbf{9.66667}	&\textbf{-125}&	\textbf{-0.54683}&	\textbf{-0.51938}\\
11.68421	&-125&	-1.15869&	-0.45924\\
3.67925	&-125	&-0.49998&	-0.5281\\
5.21754	&-125&	-0.60928&	-0.88283\\
\end{tabular}
\end{table}

\end{appendix}

\end{document}